\documentclass{amsart}
\usepackage{amssymb,amsmath}
\newtheorem{theorem}{Theorem}
\newtheorem{corollary}[theorem]{Corollary}
\newtheorem{definition}[theorem]{Definition}
\newtheorem{proposition}[theorem]{Proposition}
\newtheorem{example}[theorem]{Example}
\newtheorem{lemma}[theorem]{Lemma}
\newtheorem{remark}[theorem]{Remark}

\newcommand{\s}{\,^*\!}

\begin{document}

\baselineskip=18pt
\begin{center}
{\textbf{L\'{e}vy Processes on a First Order Model}}
\footnote{\noindent Address: School of Mathematical Sciences,
University of KwaZulu-Natal, Pietermaritzburg, 3209 South Africa;
e-mail: ngs@ukzn.ac.za }
\end{center}

\begin{center}

Siu-Ah Ng

\medskip

\emph{University of KwaZulu-Natal, South Africa}

\end{center}

\bigskip

\begin{quote}
\small{\em The classical notion of L\'{e}vy process is generalized to one that takes as its values probabilities on a first order model equipped with a commutative semigroup. This is achieved by applying a convolution product on definable probabilities and the infinite divisibility with respect to it.}

\bigskip

{\sc Key words}: First order model, independence property, L\'{e}vy process, infinite divisibility, convolution, semigroup, Loeb measure, nonstandard analysis.

{\sc Mathematics Subject Classification}: 03C99, 03H05, 60G51\\

\end{quote}

\vskip 30pt

\section{Introduction}\label{intro}

In a first order model $\mathfrak{A}=(A, \dots ),$ an element $a\in A$ is identifiable with its type in the diagram language $\mathcal{L}_A.$ As a type, $a$ corresponds to a $\{0,1\}$-probability measure on the Boolean algebra of $\mathcal{L}_A$-formulas having only one free variable. In this respect, elements in $A$ are regarded as deterministic. Then an arbitrary probability measure on the Boolean algebra corresponds to the law of a certain random variable---both the law and the random variable are liberally identified with each other. So intuitively one treats the collection of these probabilities as random ``elements" of $A$ and a time evolution of these probabilities as a stochastic process on $\mathfrak{A}.$

In classical stochastic analysis, L\'{e}vy processes are stochastic processes with stationary independent increments. Prominent examples from this important and well-studied class of processes include Brownian motions and Poisson processes. A fundamental characterization of this class is that the laws of these processes satisfy some infinitely divisible condition.

The theme of this article is to explore a way to define L\'{e}vy processes on first order models as closely analogous to the classical real-valued L\'{e}vy processes as possible. To achieve this, we need to borrow Keisler's notion of probability measure on a model and nonforking product from \cite{Kei1} and rely heavily on results from that article. However, in contrast to \cite{Kei1}, we will neither deal with forking nor delicate extensions on larger fragments. In order to specify a convolution product of two definable probabilities, we require that a commutative semigroup be definable. Then the remaining task is to identify definable probabilities which are infinitely divisible \emph{w.r.t.} the convolution product and to define L\'{e}vy processes \emph{w.r.t.} either a discrete or a continuous timeline. In either case, the dynamics behind the evolution is determined by the definable semigroup. In fact, given an infinitely divisible probability $\mu,$ the corresponding L\'{e}vy process is a time evolution from $0,$ a special element in $A,$ to $\mu.$ When $\mathfrak{A}$ is the real ordered field $(\mathbb{R},+,\cdot,\leq,0,1),$ everything here coincides with the classical real-valued L\'{e}vy processes.

Comparing with Keisler's work \cite{Kei2} on randomization of a first order model, here we only deal with objects closely connected to $\mathfrak{A}$ and will not involve an external probability space. Consequently, instead of random variables, we work purely with probabilities on $\mathfrak{A}.$ Moreover, our measure algebra already has enough saturation built into it, hence we are able to avoid technicalities such as finite additivity \emph{vs}. $\sigma$-additivity and liftings \emph{vs}. standard parts. Obviously, since we are moving away from the classical stochastic setting, a lot of analytic tools such as Fourier transforms have to be given up. One needs to find algebraic, model-theoretic and combinatorial replacements in order to obtain useful results.

Other equivalent formulations of our L\'{e}vy processes should be possible. For example, by defining hyperfinite random walks on $\mathfrak{A}$ or by starting from nonstandard compound Poisson processes. But we will not take such routes here. Interestingly, it is unclear at this point what corresponds to a Brownian motion, the prototype of L\'{e}vy processes, on a general $\mathfrak{A}.$ For further investigation, perhaps one should also study Markov processes on some first order models.

We first introduce our notation and terminologies in the next section. The role of Borel functions is played by definability in our context. The details are given in \S\ref{conv}. In \S\ref{infdiv}, infinitely divisible probabilities and L\'{e}vy processes are developed. In our context, a L\'{e}vy process can be regarded as an evolution along a ``straight line segment of probabilities" from a fixed deterministic element to a fixed infinitely divisible definable probability. The process is indexed by various types of timelines, discrete or continuous. In order to define continuous time indexed L\'{e}vy processes, convolution exponentials are introduced and the L\'{e}vy-Khintchine property is formulated. A class of Loeb measures constructed from internally definable probabilities is required and L\'{e}vy processes for these probabilities are also investigated.

\section{Basic notion and assumptions}\label{notion}

\noindent Some familiarity with model theory, stochastic analysis and nonstandard analysis is assumed. The latter is only needed for the Loeb measure construction to get $\sigma$-additive probability measures. Notation, definitions and basic results from \cite{CK}, \cite{AFHL} and \cite{Stroock} are freely used throughout.

The real closed ordered field $(\mathbb{R}, +, \cdot,\leq, 0, 1)$ is denoted by $\mathfrak{R}.$

We work with a fixed first order language $\mathcal{L}$ and an $\mathcal{L}$-model $\mathfrak{A}=(A, \dots ).$ Moreover, either we require that $\mathcal{L}$ be countable, as in \cite{Kei1}, or we we allow $\vert\mathcal{L}\vert$ be arbitrary but require ${\rm Th}(\mathfrak{A}),$ \emph{i.e.} the theory of $\mathfrak{A},$ do not have the independence property. Under the absence of the independence property, by the Corollary in \cite{Ng2}, results from \cite{Kei1} for measures on $\mathfrak{A}$ remain valid.

Furthermore, we assume that a commutative semigroup structure is definable in $\mathfrak{A}.$ That is, there is an $\mathcal{L}$-formula $\theta (x,y,z)$ such that $\mathfrak{A}$ satisfies the following:

$\left\{
\begin{array}{l}
 \hbox{$\forall x  y \exists ! z \theta (x,y,z)$}\\
 \hbox{$\forall x  y z  \big(\theta (x,y, z) \leftrightarrow \theta (y, x, z)\big)$} \\
 \hbox{$\forall x  y  z  w\Big(\big(\exists v\big(\theta (x, y, v)\land\theta (v, z, w)\big)\big)\leftrightarrow \big(\exists u\big(\theta (y, z, u)\land\theta (x, u, w)\big)\big)\Big)$}\\
 \hbox{$\exists x  \forall y\theta (x, y, y).$}
\end{array}
\right.$

For example, if $\mathfrak{A}$ defines a poset in which there is a least element and any two elements have a unique least upper bound then we can take $\theta (x, y, z)$ to be the formula saying that $z$ is the least upper bound of $x$ and $y.$

Hereafter we fix such formula $\theta.$ Write $x+y=z$ instead of $\theta (x,y,z)$ and denote the neutral element in $\mathfrak{A}$ given by the last axiom (which is necessarily unique) as $0,$ \emph{i.e.} $\mathfrak{A}\models\forall x \theta (0, x, x).$

For our purpose, we will mostly use $+$ to define an iterated convolution product of a fixed probability with itself, hence commutativity is not essential; but the notation becomes somewhat simplified and natural under this requirement.

We assume that there is an uncountable inaccessible cardinal $\kappa>\lvert\mathcal{L}\rvert+\lvert A\rvert$ and work with a $\kappa$-saturated nonstandard universe (in the sense of nonstandard analysis) containing a saturated elementary extension of $\mathfrak{A}$ of cardinality $\kappa.$ Elements in the nonstandard universe are referred to as internal objects and every standard object $X$ is extended to an internal one denoted by ${^*X}.$ So the internal model $^*\mathfrak{A}$ is the saturated elementary extension of $\mathfrak{A}$ of cardinality $\kappa.$

Sets of cardinality $<\kappa$ are called small.

Internal symbols in $^*\mathcal{L}\setminus\mathcal{L}$ will not be used.

The use of the $\kappa$ and the saturated model is purely for convenience and for consistency with the framework in \cite{Kei1}, as $\aleph_1$-saturation of the nonstandard universe would be sufficient for the Loeb measure construction and as the results needed from \cite{Kei1} can be rephrased for a $\aleph_1$-saturated elementary extension of $\mathfrak{A}.$

We mainly work with formulas in $ \mathcal{L}_A,$ \emph{i.e.} $ \mathcal{L}$ expanded by adding a new constant symbol for each element of $A,$ and regard $\mathfrak{A}$ as a $\mathcal{L}_A$-model in a canonical way. Given an $\mathcal{L}_A$-formula $\varphi=\varphi (\bar{x})$ in free variables $\bar{x}$ of arity $n,$ we let $\varphi^{A}$ denote the set $\{\bar{a}\in A^n \mid \mathfrak{A}\models \varphi(\bar{a})\}.$ (Note that $^*\varphi^{A}$ is the internal set $\{\bar{a}\in {{\s A}^n} \mid {^*\mathfrak{A}}\models \varphi(\bar{a})\},$ \emph{i.e.} $\varphi^{{\s A}}.$) Then $\mathcal{B}_0^n$ stands for the set algebra generated by the $\varphi^{A},$ where $\varphi$ ranges over $\mathcal{L}_A$-formulas of arity $n$ (for its free variables).
The notation ${\mathcal{B}^n}$ stands for the set algebra generated by the $^*\varphi^{A},$ still ranging over $\mathcal{L}_A$-formulas of arity $n.$ As Boolean algebras, ${\mathcal{B}_0^n}$ and ${\mathcal{B}^n}$ are isomorphic. Although elements in ${\mathcal{B}^n}$ are internal, ${\mathcal{B}^n}$ is in general an external subalgebra of $^*\mathcal{B}_0^n.$ The $\sigma$-algebra generated by ${\mathcal{B}^n}$ is denoted by $\sigma{\mathcal{B}^n}.$ We write $\mathcal{B}$ and $\sigma\mathcal{B}$ for $\mathcal{B}^1$ and $\sigma\mathcal{B}^1$ respectively.

Given an internal finitely additive probability measure $\mu$ on $^*\mathcal{B}_0^n,$ the Loeb measure construction from nonstandard analysis produces a $\sigma$-additive probability measure on the Loeb algebra of $^*\mathcal{B}_0^n$ (a certain $\sigma$-algebra extending the $\sigma$-algebra generated by $^*\mathcal{B}_0^n$) denoted by $L(\mu).$ Moreover, on $^*\mathcal{B}_0^n,$ values of $\mu$ and $L(\mu)$ are infinitely close to each other. $L(\mu)$ is the unique $\sigma$-additive extension \emph{w.r.t.} these properties.

Here we temporarily denote the restriction of $L(\mu)$ to $\sigma{\mathcal{B}^n}$ by ${^\circ\mu}.$

$\mathcal{M}^n$ is the notation for the $\sigma$-additive probability measures ${^\circ\mu}$ on $\sigma{\mathcal{B}^n}$ obtained in this way. $\mathcal{M}^1$ is just written as $\mathcal{M}.$ Note that if $\nu$ is a finitely additive probability measure on $\mathcal{B}_0^n$ then ${^\circ({^*\nu})}\in\mathcal{M}^n.$ Conversely, if $\mu\in\mathcal{M}^n,$ then $\mu ={^\circ({^*\nu})},$ where $\nu$ is the finitely additive probability on ${\mathcal{B}_0^n}$ given by $\nu(\phi^A)=\mu({^*\varphi}^A).$ Moreover, for $\mu\in\mathcal{M}^n$ and $X\in\sigma{\mathcal{B}^n},$ we have $\mu(X)=\sup\{\mu(Y)\mid Y\in \mathcal{B}^n\land Y\subset X\}=\inf\{\mu(Y)\mid Y\in \mathcal{B}^n\land Y\supset X\},$ that is, $\mu$ is in agreement with both the inner and the outer measure of its restriction on $\mathcal{B}^n.$ Hence elements in $\mathcal{M}$ correspond to the unique $\sigma$-additive extension of finitely additive probability measures on $\mathcal{B}^n.$ In particular, we only need to specify a finitely additive probability measure on either $\mathcal{B}_0^n$ or $\mathcal{B}^n$ in order to determine a required one in $\mathcal{M}^n.$ In general, $\mathcal{M}^n$ is a proper subclass of $\sigma$-additive probability measures on $\sigma{\mathcal{B}^n}.$

We also regard ${\s A}$ as topologized by basic open sets from $\mathcal{B}$ and when a function $f: {{\s A}^n}\longrightarrow \mathbb{R}$ is measurable \emph{w.r.t.} $\sigma{\mathcal{B}^n},$ we say that $f$ is $\sigma{\mathcal{B}^n}$-Borel. So $f$ is $\sigma{\mathcal{B}^n}$-Borel if $f^{-1}\big[(-\infty,r]\big]\in\sigma{\mathcal{B}^n}$ for every $r\in\mathbb{Q}.$

Given $\mu\in{\mathcal{M}^n}$ and $\sigma{\mathcal{B}^n}$-Borel $f$ we write $\int f(\bar{x})d\mu(\bar{x})$ for the Lebesgue integral of $f$ \emph{w.r.t.} $\mu.$

Given small $\mathfrak{C}\prec{^*\mathfrak{A}},$  all the above remain valid with $\mathfrak{A}$ replaced by $\mathfrak{C}.$ Moreover, we use notation like $\mathcal{B}_{\mathfrak{C}},\mathcal{M}_{\mathfrak{C}},\dots$ to denote the counterparts of $\mathcal{B},\mathcal{M},\dots$ for $\mathfrak{C}.$

\begin{example}
(i) For each $a\in A,$ we let $\delta_a $ denote the Dirac measure at $a.$ \emph{i.e.} for all $X\in\sigma\mathcal{B},$ $\delta_a (X)=1\,$ whenever $a\in X\,$ and $=0$ otherwise. Similarly, for $a\in{\s A},$ we let $\Delta_a$ denote the internal Dirac measure at $a,$ \emph{i.e.} for all $X\in{^*\mathcal{B}}_0,$ $\Delta_a (X)=1\,$ whenever $a\in X\,$ and $=0$ otherwise.

Note that $\delta_a\in\mathcal{M}$ for every $a\in A,$ as ${^\circ{\Delta_a}}=\delta_a.$

(ii) Now suppose $\mathcal{L} = \{\leq \}$ and $\mathfrak{A} = (\mathbb{Q},\leq ).$ Then for every infinitesimal $a\in{\s A}$ we have ${^\circ\Delta}_a =\delta_0.$\nolinebreak\hfill $\Box$
\end{example}

Probability measures in $\mathcal{M}$ are also referred to as probabilities on $\mathfrak{A}.$ We also speak of probabilities on $\mathfrak{C}$ for other $\mathfrak{C}\prec{^*\mathfrak{A}}$ in a similar manner.

Intuitively, one regards a probability on $\mathfrak{A}$ as the law of a random variable taking values in $A.$ Hence one also regard $\mathcal{M}$ as a set of random elements of $A,$ in which the deterministic ones are identified with  $\delta_a, a\in A,$ \emph{i.e.} the $\mathcal{L}_A$-types.

\section{Definable probabilities and the convolution product}\label{conv}

Let $\mathcal{A}$ be the $\sigma$-algebra generated from sets which are small union of the $\varphi^{{\s A}}$'s, where the $\varphi(x)$'s are $\mathcal{L}_{{\s A}}$-formulas. By a probability measure over $^*\mathfrak{A},$ we mean a $\sigma$-additive probability measure on $\mathcal{A}$ such that for every small set $S$ of the $\varphi^{{\s A}}$'s, there is a countable subset $S_0\subset S$ such that $\cup S$ and $\cup S_0$  have the same measure. By \cite{Kei1}~Lem.~6.1, given a probability measure $\mu$ over $^*\mathfrak{A},$ the restriction of $\mu$ to $\sigma\mathcal{B}_{\mathfrak{C}}$ belongs to $\mathcal{M}_{\mathfrak{C}}$ for any small $\mathfrak{C}\prec{^*\mathfrak{A}}.$

\begin{definition}\label{countablydefinable}
(i) Let $\mathfrak{C}\prec{^*\mathfrak{A}}$ be small. A probability measure $\mu$ over $^*\mathfrak{A}$ is called definable over $\mathfrak{C}$ if for every $\mathcal{L}$-formula $\varphi(x,\bar{y}),$ the mapping $f:{\s A}^n\longrightarrow [0,1]$ given by $f(\bar{y})=\mu\big(\varphi^{\s A}(\cdot,\bar{y})\big)$ is $\sigma\mathcal{B}_{\mathfrak{C}}^n$-Borel, where $n$ is the arity of $\bar{y}$ and $\varphi^{\s A}(\cdot,\bar{y})$ denotes $\{x\in{\s A}\mid {^*\mathfrak{A}}\models \varphi(x,\bar{y})\}.$

(Note that such $f$ is still $\sigma\mathcal{B}_{\mathfrak{C}}^n$-Borel even when $\varphi(x,\bar{y})$ is an $\mathcal{L}_C$-formula, because we can write $\varphi(x,\bar{y})$ as $\theta(x,\bar{y},\bar{c}),$ for some $\mathcal{L}$-formula $\theta$ and some $\bar{c}$ from $C$ of arity $m,$ then $f$ is the section of a $\sigma\mathcal{B}_{\mathfrak{C}}^{n+m}$-Borel function at $\bar{c},$ hence is $\sigma\mathcal{B}_{\mathfrak{C}}^n$-Borel.)

(ii) The probability measure $\mu$ over $^*\mathfrak{A}$ is said to be countably definable over $\mathfrak{C}$ if, for each $n\in\mathbb{N},$ there is a countably generated subalgebra $\mathcal{C}_n\subset\sigma\mathcal{B}_{\mathfrak{C}}^n,$ such that the mapping $\bar{y}\longmapsto \mu\big(\varphi^{\s A}(\cdot,\bar{y})\big)$ is $\mathcal{C}_n$-Borel for any $\mathcal{L}$-formula $\varphi(x,\bar{y})$ of arity $n.$ \nolinebreak\hfill $\Box$
\end{definition}
We obtain from \cite{Kei1}~Prop.~6.4(ii) and Cor.~6.7 the following.

\begin{lemma}\label{definable1} Suppose ${\rm Th}(\mathfrak{A})$ does not have the independence property. Then every $\mu\in\mathcal{M}$ extends to a countably definable probability measure $\bar{\mu}$ over $^*\mathfrak{A}.$\nolinebreak\hfill $\Box$
\end{lemma}

In general, the extension $\bar{\mu}$ above is not unique. For example, if there is a type over $\mathcal{L}_A$ omitted by $\mathfrak{A}$ but realized by two distinct $a,b\in{\s A},$ then, considering the Dirac measures on $\mathcal{A},$ both $\delta_a$ and $(\delta_a+\delta_b)/2$ are countably definable probability measures over $^*\mathfrak{A}$ extending their common restriction on $\mathfrak{A}.$

However if the extension is definable over $\mathfrak{A}$ then there is uniqueness in the following sense for $\mathcal{L}_A$-formulas.

\begin{lemma}\label{definable2}
Let $\mu_1,\mu_2$ be probability measures over $^*\mathfrak{A}$ extending some $\mu\in\mathcal{M}.$ Suppose both $\mu_1$ and $\mu_2$ are definable over $\mathfrak{A}.$ Then for each $\nu\in\mathcal{M}^n$ and $\mathcal{L}_A$-formula $\varphi(x,\bar{y}),$ where $n$ is the arity of $\bar{y},$
\[\nu\big(\{\bar{y}\in{\s A}^n\mid \mu_1(X_{\bar{y}})=\mu_2(X_{\bar{y}})\}\big)=1,\quad\text{where}\;X_{\bar{y}}=\varphi^{\s A}(\cdot,\bar{y}).\]
\end{lemma}

\begin{proof}
Suppose on the contrary, $\nu\big(\{\bar{y}\in{\s A}^n\mid \mu_1(X_{\bar{y}})\neq\mu_2(X_{\bar{y}})\}\big)>0.$ By Def.~\ref{countablydefinable} and the assumption, $\{\bar{y}\in{\s A}^n\mid \mu_1(X_{\bar{y}})\neq\mu_2(X_{\bar{y}})\}\in\sigma\mathcal{B}^n.$ So, as $\nu\in\mathcal{M}^n,$ there is $\mathcal{L}_A$-formula $\theta(\bar{y})$ such that $\nu(\theta^{{\s A}})>0$ and $\theta^{\s A}\subset\{\bar{y}\in{\s A}^n\mid \mu_1(X_{\bar{y}})\neq\mu_2(X_{\bar{y}})\}.$ In particular, ${^*\mathfrak{A}}\models\exists \bar{y}\theta(\bar{y}),$ hence, as $\mathfrak{A}\prec {^*\mathfrak{A}},$ there is $\bar{a}\in A^n\cap\theta^{\s A}.$  Then $\mu_1\big(\varphi^{\s A}(\cdot,\bar{a})\big)\neq\mu_2\big(\varphi^{\s A}(\cdot,\bar{a})\big),$ contradicting to $\mu_1,\mu_2$ extending $\mu.$
\end{proof}

The following terminology deviates a bit from that used in \cite{Kei1} and it does not apply to probability measures over $^*\mathfrak{A}.$

\begin{definition}\label{definableextension}
Suppose $\mathfrak{C}\prec{^*\mathfrak{A}}$ and $\mathfrak{C}$ is small.

(i) A definable probability on $\mathfrak{C}$ is some $\mu\in\mathcal{M}_{\mathfrak{C}}$ having an extension to a probability measure $\bar{\mu}$ over $^*\frak{A}$ such that $\bar{\mu}$ is countably definable over $\mathfrak{C}.$ A probability on $\mathfrak{C}$ is simply called definable if the reference to $\mathfrak{C}$ is clear.

(ii) If the $\bar{\mu}$ is unique, we say that $\mu$ is strongly definable.

(iii) The set of definable probabilities on $\mathfrak{C}$ is denoted by $\mathcal{D}(\mathfrak{C}).$\nolinebreak\hfill $\Box$
\end{definition}

Clearly we have the following.

\begin{proposition}\label{delta0}
$\delta_a$ is strongly definable, for each $a\in A.$\nolinebreak\hfill $\Box$
\end{proposition}

\begin{remark}\label{delta1}
From Prop.~\ref{delta0}, we see that $A$ embeds into $\mathcal{D}(\mathfrak{A})$ in a canonical way via $a\longmapsto \delta_a.$ So we can view $\mathcal{D}(\mathfrak{A})$ as an expansion of $\mathfrak{A}$ by including ``definable random elements"---a sort of ``definable randomization" of $\mathfrak{A}$ and view $A$ as the set of deterministic elements in it. Moreover, the satisfaction relation $\mathfrak{A}\models \varphi(a)$ becomes $\delta_a(\varphi^{\s A})=1.$\nolinebreak\hfill $\Box$
\end{remark}

In \cite{Kei1} a class of probability measures called ``smooth measures" was studied. These are probability measures that exclude distinct extensions on the unstable part, thus generalize stable types in the classical theory. By \cite{Kei1}~Prop.~6.4(iii) and and Cor.~6.7, we have the following.

\begin{lemma}\label{definable3}
{\rm (i)} Each smooth probability on $\mathfrak{A}$ is strongly definable.

{\rm (ii)} In particular, if ${\rm Th}(\mathfrak{A})$ is stable, all probabilities on $\mathfrak{A}$ are strongly definable, so $\mathcal{M}=\mathcal{D}(\mathfrak{A}).$ \nolinebreak\hfill $\Box$
\end{lemma}

In the absence of the independence property, every probability on $\mathfrak{A}$ extends to a smooth one by \cite{Kei1}~Thm.~3.16(ii). Therefore the following holds.

\begin{lemma}\label{definable4} Suppose ${\rm Th}(\mathfrak{A})$ does not have the independence property.
Let $\mu\in\mathcal{M},$ then for some small $\mathfrak{C}$ with $\mathfrak{A}\prec\mathfrak{C}\prec{^*\mathfrak{A}},$  $\mu$ has a extension to a strongly definable probability on $\mathfrak{C}.$  \nolinebreak\hfill $\Box$
\end{lemma}

If $\mathfrak{A}$ is o-minimal, then its theory does not have the independence property. Hence the above lemmas are applicable to o-minimal models. In these models, a certain linear order is definable and every $\mathcal{L}_A$-formula $\varphi (x)$ is equivalent to a finite combination of intervals. Important examples include $ \mathfrak{R}$ and its expansions equipped with the exponential function or restricted analytic functions. In particular, the law of a real-valued random variable in the classical sense is always strongly definable. Therefore our setting here is a generalization of that for classical stochastic analysis.

It is worth mentioning that $p$-adic fields are not o-minimal but does not have the independence property either.

But we do not know whether the absence of the independence property or o-minimality or elimination of quantifiers imply that every $\mu\in\mathcal{M}$ is always strongly definable, or even just definable.

A useful fact is the following that $\mathcal{D}(\mathfrak{A})$ and its subclass of strongly definable probabilities are closed under convex combinations. The verification is straight-forward.

\begin{proposition}\label{lincon}
Let $\mu$ and $\nu$ be (strongly) definable probabilities on $\mathfrak{A}.$ Let $r\in[0,1].$ Then the probability $r\mu+(1-r)\nu$ is also a (strongly) definable probability on $\mathfrak{A}.$ \nolinebreak\hfill $ \Box$
\end{proposition}

\begin{theorem}\label{convolution0}
Let $\mu,\nu$ be probability measures over $^*\frak{A}$ so that both $\mu$ and $\nu$ are countably definable over some small $\mathfrak{C}\prec{^*\mathfrak{A}.}$ Then the mapping
\[\xi:\varphi^{\s A}\longmapsto\int\mu\big(\{x\in{\s A}\mid {^*\mathfrak{A}}\models\varphi(x+y)\}\big)d\nu(y),\quad\text{where}\;\varphi(x)\;\text{is an}\;\mathcal{L}_{\s A}\text{-formula,}\]
extends uniquely to a probability measure over $^*\mathfrak{A}$ which is countably definable over $\mathfrak{C}.$

We denote this measure by $\mu\bigstar\nu.$
\end{theorem}

\begin{proof}
First, by the note in Def.~\ref{countablydefinable}, the function $y\longmapsto \mu\big(\{x\in{\s A}\mid {^*\mathfrak{A}}\models\varphi(x+y)\}\big)$ is $\mathcal{A}$-Borel, hence the above Lebesgue integral, therefore the $\xi,$ is well-defined.

Clearly $\xi$ forms a finitely additive probability measure on $\mathcal{A}.$ For example, Let $\varphi_1(x),\varphi_2(x)$ be $\mathcal{L}_{\s A}$-formulas such that $\varphi_1^{\s A}\cap\varphi_2^{\s A}=\emptyset,$ then
\[\xi\big(\varphi_1^{\s A}\cup\varphi_2^{\s A}\big)=\int\mu\big((\varphi_1\vee\varphi_2)^{\s A}(\cdot +y)\big)d\nu(y)=\int\Big(\mu\big(\varphi_1^{\s A}(\cdot +y)\big)+\mu\big(\varphi_2^{\s A}(\cdot +y)\big)\Big)d\nu(y)=\xi\big(\varphi_1^{\s A}\big)+\xi\big(\varphi_2^{\s A}\big).\]
So, by \cite{Kei1}~Lem.~6.2, $\xi$ extends uniquely to probability measure $\xi^\prime$ over $\s\mathfrak{A}.$

For any $\mathcal{L}_{\s A}$-formula $\varphi(x),$ we have $\xi\big(\varphi^{\s A}\big)=[\mu\times \nu]\big(\{(x,y)\in{\s A}^2\mid {^*\mathfrak{A}}\models\varphi(x+y)\}\big),$ according to the definition of the nonforking product $[\mu\times\nu]$ given in \cite{Kei1}~Def.~6.11. Hence, by \cite{Kei1}~Lem.~6.13, $\xi^\prime$ is countably definable over $\mathfrak{C}.$
\end{proof}

For definable probabilities, the probability measure given above is unique on $\sigma\mathcal{B}$ in the following sense.

\begin{proposition}\label{uniqueconv}
Let $\mu,\nu\in\mathcal{D}(\mathfrak{A}).$ Let $\mu$ and $\nu$ respectively extend to probability measures $\mu_1,\mu_2$ and $\nu_1,\nu_2$ over $^*\mathfrak{A}$ which are countably definable over $\mathfrak{A}.$

Then $\mu_1\bigstar\nu_1$ and $\mu_2\bigstar\nu_2$ (as given by Thm.~\ref{convolution0}) coincide on $\sigma\mathcal{B}.$
\end{proposition}

\begin{proof}
Let $\varphi(x)$ be any $\mathcal{L}_A$-formula. Then
\[\int\mu_2\big(\varphi^{\s A}(\cdot +y)\big)d\nu_2(y)=\int\mu_2\big(\varphi^{\s A}(\cdot +y)\big)d\nu(y)=\int\mu_1\big(\varphi^{\s A}(\cdot +y)\big)d\nu(y)=\int\mu_1\big(\varphi^{\s A}(\cdot +y)\big)d\nu_1(y),\]
where the first equality follows from the mapping $y\longmapsto \mu_2\big(\varphi^{\s A}(\cdot +y)\big)$ being $\sigma\mathcal{B}$-Borel and $\nu_2$ extending $\nu,$ the second one follows from Lem.~\ref{definable2} and the third one follows from $\nu_1$ extending $\nu.$

Hence $\mu_1\bigstar\nu_1$ and $\mu_2\bigstar\nu_2$ coincide on $\sigma\mathcal{B}.$
\end{proof}

Now we define the convolution of definable probabilities on $\mathfrak{A}.$

\begin{definition}\label{convolution}
(i) Let $\mu,\nu\in\mathcal{D}(\mathfrak{A}).$ Let $\bar{\mu},\bar{\nu}$ be any extensions of $\mu,\nu$ respectively which are probability measures over $^*\mathfrak{A}$ and countably definable over $\mathfrak{A}.$ Then the convolution product of $\mu$ and $\nu$ is the restriction of the $\bar{\mu}\bigstar\bar{\nu}$ given by Thm.~\ref{convolution0} to $\sigma\mathcal{B}.$

(ii) By Prop.~\ref{uniqueconv} this convolution product is unique. This uniquely defined probability measure on $\mathfrak{A}$ is denoted by $\mu\star\nu.$\nolinebreak\hfill $\Box$
\end{definition}

By Thm.~\ref{convolution0}, we also have the following.

\begin{corollary}\label{convolution2}
Let $\mu,\nu\in\mathcal{D}(\mathfrak{A}).$

{\rm (i)} $(\mu\star\nu)\in\mathcal{D}(\mathfrak{A}).$

{\rm (ii)} Moreover, if $\mu,\nu$ respectively extend to probability measures $\bar{\mu},\bar{\nu}$ over $^*\mathfrak{A}$ which are countably definable over $\mathfrak{A},$ then $\mu\star\nu$ extends to $\bar{\mu}\bigstar\bar{\nu},$ a probability measures over $^*\mathfrak{A}$ which is countably definable over $\mathfrak{A}.$\nolinebreak\hfill $\Box$
\end{corollary}

\begin{remark}\label{defrmk}
(i) If the above $\mu,\nu$ are strongly definable, is $\mu\star\nu$ strongly definable?

(ii) If ${\rm Th}(\mathfrak{A})$ does not have the independence property, then by Lem.~\ref{definable4} and $\kappa>\vert\mathcal{L}\vert+\vert A\vert$ being an uncountable inaccessible cardinal, there is a small $\mathfrak{C}$ such that $\mathfrak{A}\prec\mathfrak{C}\prec{^*\mathfrak{A}}$ and every probability on $\mathfrak{A}$ extends to a strongly definable one on $\mathfrak{C}.$ By restricting the convolution product back to $\sigma\mathcal{B}$ one can define a convolution product for the whole $\mathcal{M},$ although this need not be the only possible one.

(iii) If one of the probability measures $\mu,\nu$ on $\mathfrak{A}$ extends to a smooth probability measure over $^*\mathfrak{A},$ then $\mu\star\nu=\nu\star\mu$ holds by \cite{Kei1}~Thm.~6.15 and the commutativity of $+.$

(iv) Given $\mu\in\mathcal{D}(\mathfrak{A})$ and $a\in A,$ we have $(\mu\star\delta_a)(\varphi^{\s A}(\cdot))=\mu(\varphi^{\s A}(\cdot+a)),$ \emph{i.e.} $\mu\star\delta_a$ is the probability measure given by the translation of $\mu$ by $a.$ In particular, $\mu\star\delta_0=\mu.$

(v) By using the Dirac measure at $a$ on $\mathcal{A},$ $\delta_a$ extends to a probability measure over $^*\mathfrak{A}$ which is countably definable over $\mathfrak{A}$ and is smooth, so it follows from (iii) that $\mu\star\delta_a=\delta_a\star\mu.$\nolinebreak\hfill $\Box$
\end{remark}

Since our attention is more on the iterated convolution of a single definable probability on $\mathfrak{A}$ we are not so concerned with defining a convolution product for the whole $\mathcal{M}$ nor with commutativity. However, in our case we do have the convenience of associativity of the convolution product.

\begin{proposition}\label{associative}
Let $\mu,\nu,\lambda\in\mathcal{D}(\mathfrak{A}).$ Then $(\mu\star\nu)\star\lambda =\mu\star(\nu\star\lambda).$
\end{proposition}

\begin{proof}
Let $\mu,\nu,\lambda$ respectively extend to $\bar{\mu},\bar{\nu},\bar{\lambda},$ some probability measures over $^*\mathfrak{A}$ which are countably definable over $\mathfrak{A}.$

We first note that for any $\mathcal{L}_{\s A}$-formula $\theta(x),$ $(\bar{\nu}\bigstar\bar{\lambda})\big(\theta^{\s A}(\cdot)\big)=[\bar{\nu}\times\bar{\lambda}]\big(\theta^{\s A}(\cdot+\cdot)\big),$ where $[\bar{\nu}\times\bar{\lambda}]$ denotes the nonforking product in \cite{Kei1}.

Therefore, for any $\sigma\mathcal{B}$-Borel function $f:{\s A}^{1+n}\to\mathbb{R}$ and $a_1,\dots,a_n\in{\s A},$ we have
\[\int f(w,a_1,\dots,a_n)d(\bar{\nu}\bigstar\bar{\lambda})(w)=\int f(y+z,a_1,\dots,a_n)d[\bar{\nu}\times\bar{\lambda}](y,z).\]
So
\[\int f(w,a_1,\dots,a_n)d(\bar{\nu}\bigstar\bar{\lambda})(w)=\iint f(y+z,a_1,\dots,a_n)d\bar{\nu}(y)d\bar{\lambda}(z),\]
by \cite{Kei1}~Lem.~6.13.

Now let $\varphi(x)$ be an $\mathcal{L}_A$-formula. By Cor.~\ref{convolution2}, $\bar{\nu}\bigstar\bar{\lambda}$ is a probability measure over $^*\mathfrak{A},$ countably definable over $\mathfrak{A}$ and extending $\nu\star\lambda.$ Similarly, $\bar{\mu}\bigstar\bar{\nu}$ is a probability measure over $^*\mathfrak{A},$ countably definable over $\mathfrak{A}$ and extending $\mu\star\nu.$ So from what was just shown, we have:
\begin{align*}
\big(\mu\star(\nu\star\lambda)\big)\big(\varphi^{\s A}(\cdot)\big) &=\int \bar{\mu}\big(\varphi^{\s A}(\cdot+w)\big)d(\bar{\nu}\bigstar\bar{\lambda})(w)=\iint \bar{\mu}\big(\varphi^{\s A}(\cdot+(y+z))\big)d\bar{\nu}(y)d\bar{\lambda}(z)\\
            &=\iint \bar{\mu}\big(\varphi^{\s A}((\cdot+y)+z)\big)d\bar{\nu}(y)d\bar{\lambda}(z)=\int (\bar{\mu}\bigstar\bar{\nu})\big(\varphi^{\s A}(\cdot+z)\big)d\bar{\lambda}(z)\\
            &=\big((\bar{\mu}\bigstar\bar{\nu})\bigstar\bar{\lambda}\big)\big(\varphi^{\s A}(\cdot)\big)=\big((\mu\star\nu)\star\lambda\big)\big(\varphi^{\s A}(\cdot)\big).
\end{align*}
Therefore $(\mu\star\nu)\star\lambda =\mu\star(\nu\star\lambda).$
\end{proof}

Notice the use of the associativity of $+$ in the above proof. The following useful fact is easy to check.

\begin{proposition}\label{combination}
Let $\mu,\nu,\lambda\in\mathcal{D}(\mathfrak{A})$ and $r\in[0,1].$ Then $\big(r\mu+(1-r)\nu\big)\star\lambda =r(\mu\star\lambda)+(1-r)(\nu\star\lambda).$\nolinebreak\hfill $\Box$
\end{proposition}

\begin{definition}\label{definfdiv}
Let $\mu\in\mathcal{D}(\mathfrak{A}).$ For $n\in\mathbb{N}^+\,$ we write $\mu^{n\star}$ for the convolution power $\underbrace{\mu\star\cdots\star\mu}_{n\;\text{\rm times}}.$
When $n=0,$ $\mu^{n\star}$ is defined to be $\,\delta_0.\,$ Note that, by Prop.~\ref{associative}, $\mu^{n\star}$ is unambiguously defined.

If $\mu$ is a countably definable probability measure over $^*\mathfrak{A},$ we define $\mu^{n\bigstar}$ similarly.\nolinebreak\hfill $\Box$
\end{definition}

The following generalizes part of Prop.~\ref{lincon}, by showing that $\mathcal{D}(\mathfrak{A})$ is closed under infinite convex combinations. As a consequence, $\mathcal{D}(\mathfrak{A})$ supports some functional calculus, \emph{i.e.} some ``analytic" functions are definable.

\begin{theorem}\label{expconvolution1}
Suppose $a_n\geq 0, n\in\mathbb{N},$ such that $\sum_{n=0}^\infty a_n =1.$

Let $\mu_n\in\mathcal{D}(\mathfrak{A}), n\in\mathbb{N}.$ Then $\Big(\sum_{n=0}^\infty a_n \mu_n\Big)\in\mathcal{D}(\mathfrak{A}).$

In particular, $\Big(\sum_{n=0}^\infty a_n \mu^{n\star}\Big)\in\mathcal{D}(\mathfrak{A})$ for any $\mu\in\mathcal{D}(\mathfrak{A}).$
\end{theorem}

\begin{proof}
First, it is clear that $\Big(\sum_{n=0}^\infty a_n \mu_n\Big)\in\mathcal{M}.$

For each $n\in\mathbb{N},$ let $\mu_n$ extend to some $\bar{\mu}_n$ which is a probability measure over $^*\frak{A}$ and is countably definable over $\mathfrak{A}.$ It is easy to see that $\sum_{n=0}^\infty a_n \bar{\mu}_n$ is a probability measure over $^*\mathfrak{A}.$

Let $\bar{y}$ have arity $m.$ By definition, since the $\bar{\mu}_n$'s are countably definable over $\mathfrak{A},$  there is a countably generated subalgebra $\mathcal{C}_m\subset\sigma\mathcal{B}^m$ such that for all $\mathcal{L}$-formula $\varphi(x,\bar{y})$ and $n\in\mathbb{N},$ the mappings $\bar{y}\longmapsto\bar{\mu}_n\big(\varphi^{\s A}(\cdot,\bar{y})\big)$ are $\mathcal{C}_m$-Borel. Then, as
\[\Big(\sum_{n=0}^\infty a_n \bar{\mu}_n\Big)\big(\varphi^{\s A}(\cdot,\bar{y})\big) =\sup_{m\in\mathbb{N}}\Big(\sum_{n=0}^m a_n\bar{\mu}_n\Big)\big(\varphi^{\s A}(\cdot,\bar{y})\big),\]
$\bar{y}\longmapsto\sum_{n=0}^\infty a_n\bar{\mu}_n\big(\varphi^{\s A}(\cdot,\bar{y})\big)$ is the supremum of a countable family of $\mathcal{C}_m$-Borel functions.

Hence $\bar{y}\longmapsto\sum_{n=0}^\infty a_n\bar{\mu}_n\big(\varphi^{\s A}(\cdot,\bar{y})\big)$ is a $\mathcal{C}_m$-Borel function.

Therefore $\sum_{n=0}^\infty a_n \bar{\mu}_n$ is a probability measure over $^*\frak{A}$ which is countably definable over $\mathfrak{A}.$

Since $\sum_{n=0}^\infty a_n \bar{\mu}_n$ extends $\sum_{n=0}^\infty a_n \mu_n,$ we conclude that $\Big(\sum_{n=0}^\infty a_n \mu_n\Big)\in\mathcal{D}(\mathfrak{A}).$
\end{proof}

\begin{remark}\label{integral}
When dealing with an infinite convex combination of iterated convolution powers of some $\mu\in\mathcal{D}(\mathfrak{A}),$ the Borel functions in the above proof can be represented by integrals. Let $\bar{\mu}$ be an extension of $\mu$ which is a probability measure over $^*\frak{A}$ and is countably definable over $\mathfrak{A}.$ Then
\begin{align*}
\Big(\sum_{n=0}^\infty a_n \bar{\mu}^{n\bigstar}\Big)\big(\varphi^{\s A}(\cdot,\bar{y})\big) &=\sum_{n=0}^\infty a_n\bar{\mu}^{n\bigstar}\big(\varphi^{\s A}(\cdot,\bar{y})\big)\\
   &=\sum_{n=1}^\infty a_n\int\bar{\mu}^{(n-1)\bigstar}\big(\varphi^{\s A}(\cdot+x,\bar{y})\big) d\bar{\mu}(x)+a_0\delta_0\big(\varphi^{\s A}(\cdot,\bar{y})\big)\\
   &=\int\Big[\sum_{n=1}^\infty a_n\bar{\mu}^{(n-1)\bigstar}\big(\varphi^{\s A}(\cdot+x,\bar{y})\big)\Big] d\bar{\mu}(x)+a_0\delta_0\big(\varphi^{\s A}(\cdot,\bar{y})\big),
\end{align*}
by the Lebesgue monotone convergence theorem.\nolinebreak\hfill $\Box$
\end{remark}

We now introduce the notion of convolution exponential, in order to obtain infinitely divisible definable probabilities on $\mathfrak{A}$ in the next section.

\begin{definition}\label{expconvolution}
Let $\mu\in\mathcal{D}(\mathfrak{A}), r\in [0,\infty)$ then we define the convolution exponential of $r\mu $ as
\[e^{r(\mu\star -1)} = e^{-r}\sum_{n=0}^\infty\frac{r^n}{n!}\mu^{n\star}.\]
(So it gives $\delta_0$ when $r=0.$)

If $\mu$ is a countably definable probability measure over $^*\mathfrak{A},$  we define $e^{r(\mu\bigstar -1)}$ similarly. \nolinebreak\hfill $ \Box$
\end{definition}

As a consequence of Thm.~\ref{expconvolution1}, we have the following.

\begin{corollary}
Let $\mu\in\mathcal{D}(\mathfrak{A}), r\in [0,\infty),$ then $e^{r(\mu\star -1)}\in\mathcal{D}(\mathfrak{A}).$

Similarly, if $\mu$ is a countably definable probability measure over $^*\mathfrak{A},$  $e^{r(\mu\bigstar -1)}$ is a countably definable probability measure over $^*\mathfrak{A}.$\nolinebreak\hfill $\Box$
\end{corollary}

\section{Infinitely divisible probabilities and L\'{e}vy processes}\label{infdiv}

In the this section, we will work with probabilities from $\mathcal{D}(\mathfrak{A})$ and the convolution product $\star$ defined on it. By Cor.~\ref{convolution2}, $\mathcal{D}(\mathfrak{A})$ it is closed under $\star$ and, by Thm.~\ref{expconvolution1}, it is closed under infinite convex combinations, including the convolution exponentials. Moreover, $\mathcal{D}(\mathfrak{A})$ includes $A$ in the sense that $\delta_a\in\mathcal{D}(\mathfrak{A})$ for every $a\in A.$

From now on, $\mathcal{D}(\mathfrak{A}),$ $\star,$ as well as the elements of the underlying $\sigma\mathcal{B}$ are considered to be standard objects and, for the purpose of Loeb measure construction, we shall adopt an $\aleph_1$-saturated nonstandard universe possibly different from the one used in the previous sections. However, we continue to use $^*X$ to denote nonstandard extensions in this nonstandard universe. In particular, $^*\mathfrak{A}$ refers to the $\aleph_1$-saturated elementary extension of $\mathfrak{A}$ in this nonstandard universe, not the one used in previous sections.

Let $\nu$ be an internal probability on $^*\mathfrak{A}.$ In particular, $\nu$ is an internal finitely additive probability measure on $^*\sigma\mathcal{B}.$ So its Loeb measure $L(\nu)$ is a $\sigma$-additive probability measure on $L({^*\sigma\mathcal{B}}),$ the Loeb algebra of $^*\sigma\mathcal{B}.$ Note that $\sigma\mathcal{B}$ embeds canonically as a subalgebra of $L({^*\sigma\mathcal{B}})$ and we let $^\circ\nu$ denote the restriction of $L(\nu)$ to $\sigma\mathcal{B}.$ Note also that for $\mu\in\mathcal{M},$ $^\circ({^*\mu})=\mu.$ But, in general, $^\circ\nu$ needs not be an element of $\mathcal{M}.$

It seems necessary to expand $\mathcal{D}(\mathfrak{A})$ to a larger class of probability measures in order to make more random elements available for constructing L\'{e}vy processes and for simplifying the proof of some properties about $\mathcal{D}(\mathfrak{A})$ itself.

\begin{definition}\label{neardefinable}
(i) For $\nu_1,\nu_2\in{^*\mathcal{D}}(\mathfrak{A}),$ we write $\nu_1\approx\nu_2$ if $\nu_1(X)\approx\nu_2(X)$ for all $X\in{^*\sigma\mathcal{B}}.$ We also use expressions like $\nu_1\leq\nu_2$ and $\nu_1\lessapprox\nu_2$ in a similar manner.

(ii) A $\sigma$-additive probability measure $\mu$ on $L({^*\sigma\mathcal{B}})$ is called $L$-definable if $\mu=L(\nu)$ for some $\nu\in{^*\mathcal{D}}(\mathfrak{A}).$

(iii) We write $\Gamma({^*\mathfrak{A}})=\{L(\nu)\mid\nu\in{^*\mathcal{D}}(\mathfrak{A})\},$ \emph{i.e.} the $L$-definable probability measures on $L({^*\sigma\mathcal{B}}).$

(iv) On $\Gamma({^*\mathfrak{A}}),$ the convolution product is defined by $L(\nu_1)\star L(\nu_2)=L(\nu_1\star\nu_2).$

(Note the multiple usages of the symbol $\star:$ on ${\mathcal{D}}(\mathfrak{A}), {^*\mathcal{D}}(\mathfrak{A})$ or $\Gamma({^*\mathfrak{A}}),$ depending on the context.) \nolinebreak\hfill $ \Box$
\end{definition}

Observe that each $\mu\in\mathcal{D}(\mathfrak{A})$ extends uniquely to $L(^*\mu)\in\Gamma({^*\mathfrak{A}}).$

\begin{proposition}\label{neardefinable1}
{\rm (i)} The convolution product given in Def.~\ref{neardefinable}(iv) is well-defined.

{\rm (ii)} Given a series of the form in Thm.~\ref{expconvolution1}, we have $L\Big(\sum_{n=0}^\infty a_n \mu_n\Big)=\sum_{n=0}^\infty a_n L(\mu_n)$ for any $\mu_n\in{^*\mathcal{D}}(\mathfrak{A}), n\in\mathbb{N}.$ In particular, $L\big(e^{r(\mu\star -1)}\big)=e^{r(L(\mu)\star -1)}$ for any $\mu\in{^*\mathcal{D}}(\mathfrak{A})$ and $r\in [0,\infty).$

{\rm (iii)} $\Gamma({^*\mathfrak{A}})$ is closed under the convolution product and infinite convex combinations as in Thm.~\ref{expconvolution1}. (Hence is closed under convolution exponentials.) Moreover, ${^*\mathfrak{A}}$ embeds into $\Gamma({^*\mathfrak{A}})$ via $a\longmapsto\delta_a,$ where $a\in{\s A}.$
\end{proposition}

\begin{proof}
{\rm (i)}: Let $\mu_1,\mu_2,\nu_1,\nu_2\in{^*\mathcal{D}}(\mathfrak{A}),$ such that $L(\mu_1)=L(\nu_1)$ and $L(\mu_2)=L(\nu_2).$ Then by the Loeb construction, $\mu_1\approx\nu_1$ and $\mu_2\approx\nu_2.$

Let $X\in{^*\sigma\mathcal{B}}.$ Then by $\mu_1,\mu_2,\nu_1,\nu_2\in{^*\mathcal{M}},$ some $Y\in {^*\mathcal{B}}$ can be chosen such that
\[(\mu_1\star\mu_2)(X)\approx (\mu_1\star\mu_2)(Y)\quad\text{and}\quad (\nu_1\star\nu_2)(X)\approx (\nu_1\star\nu_2)(Y).\]
By the transfer principle from nonstandard analysis, $Y$ is represented by some ${^*{\mathcal{L}}}_{A}$-formula $\varphi$ as the set $\varphi^{C}(\cdot)$ for some internal $\mathfrak{C}\succ\mathfrak{A}.$ (In fact $\mathfrak{C}$ has the form ${^*{\overline{\mathfrak{A}}}}$ if we let ${\overline{\mathfrak{A}}}$ denote the $^*\mathfrak{A}$ used in previous sections.) Then
\begin{align*}
(\mu_1\star\mu_2)(X)&\approx (\mu_1\star\mu_2)\big(\varphi^{C}(\cdot)\big)=\int\mu_1\big(\varphi^{C}(\cdot+x)\big)d\mu_2(x)\\
                      &\approx\int {\rm st}\Big[\mu_1\big(\varphi^{C}(\cdot+x)\big)\Big]dL(\mu_2)(x)=\int L(\mu_1)\big(\varphi^{C}(\cdot+x)\big)dL(\mu_2)(x),
\end{align*}
because the integrand $\mu_1\big(\varphi^{C}(\cdot+x)\big)$ is $S$-integrable, so it lifts the Loeb integrable function ${\rm st}\Big[\mu_1\big(\varphi^{C}(\cdot+x)\big)\Big].$ (Here ${\rm st}$ denotes the standard part of a hyperreal number.) The last equality follows from the Loeb measure construction.

Likewise, $(\nu_1\star\nu_2)(X)\approx\int L(\nu_1)\big(\varphi^{C}(\cdot+x)\big)dL(\nu_2)(x),$ hence $(\mu_1\star\mu_2)(X)\approx(\nu_1\star\nu_2)(X).$ Since this holds for all $X\in{^*\sigma\mathcal{B}},$ we have $(\mu_1\star\mu_2)\approx(\nu_1\star\nu_2)$ and conclude that $L(\mu_1\star\mu_2)=L(\nu_1\star\nu_2).$

{\rm (ii)}: By a proof similar to that of Thm.~\ref{expconvolution1}. Note that $\sum_{n=0}^\infty a_n \mu_n$ refers to an internal series extending the given one.

{\rm (iii)} is clear from the corresponding properties of ${^*\mathcal{D}}(\mathfrak{A})$ and the Loeb measure construction and (ii).
\end{proof}

\begin{remark}\label{neardefrmk}
In a way, the relation between ${^*\mathfrak{A}}$ and $\Gamma({^*\mathfrak{A}})$ is like that between $\mathfrak{A}$ and ${\mathcal{D}}(\mathfrak{A}),$ so we should view $\Gamma({^*\mathfrak{A}})$ as the set of definable random elements from ${^*\mathfrak{A}}.$\nolinebreak\hfill $\Box$
\end{remark}

We will study the infinite divisibility of a probability from $\mathcal{D}(\mathfrak{A})$ or from $\Gamma({^*\mathfrak{A}})$ and the L\'{e}vy process corresponding to such probability. As remarked before, this setting generalizes the setting for classical stochastic analysis, which corresponds to the case $\mathfrak{A}=\mathfrak{R}.$ Results about classical stochastic analysis can be found in \cite{BMR}, \cite{Sato} and \cite{Stroock}, while nonstandard treatment of them can be found in \cite{AH}, \cite{Lind} and \cite{Ng2}.

\begin{definition}\label{definfdiv2}
Let $\mu\in\mathcal{D}(\mathfrak{A}).$ $\mu$ is said to be infinitely divisible if for every $n\in\mathbb{N}$ there is $\nu\in\mathcal{D}(\mathfrak{A})$ such that $\nu^{n\star} =\mu.$

This is also similarly defined for $\mu\in\Gamma({^*\mathfrak{A}}).$

To emphasize, we say ``infinitely divisible in $\mathcal{D}(\mathfrak{A})$" or ``infinitely divisible in $\Gamma({^*\mathfrak{A}})$".\nolinebreak\hfill $\Box$
\end{definition}

In the case $\mathfrak{A} = \mathfrak{R},$ classical result shows that the above $\nu$ is unique if it exists. However the uniqueness depends at least on the underlying semigroup structure of $\mathfrak{A}.$ For example, if the underlying semigroup is $[0,1)$ with $x+y$ defined to be the fractional part of the usual addition, then for $\mu=\delta_0$ the above $\nu$ would not be unique.

\begin{proposition}\label{infdiv1}
The following are equivalent for any $\mu\in{^*\mathcal{D}}(\mathfrak{A}):$

{\rm (i)} $L(\mu)$ is infinitely divisible in $\Gamma({^*\mathfrak{A}}).$

{\rm (ii)} For all small infinite $N\in{^*\mathbb{N}},$ there is $\nu\in{^*\mathcal{D}(\mathfrak{A})}$ such that $\mu \approx\nu^{N\star}.$

{\rm (iii)} For some infinite $N\in{^*\mathbb{N}},$ there is $\nu\in{^*\mathcal{D}(\mathfrak{A})}$ such that $\mu \approx\nu^{N!\star}.$
\end{proposition}

\begin{proof}
(i) $\Rightarrow$ (ii): For each $n\in\mathbb{N},$ let $\mu_n\in{^*\mathcal{D}}(\mathfrak{A})$ such that $L(\mu)=L({\mu_n})^{n\star}.$ By Def.~\ref{neardefinable}(iv) and Prop.~\ref{neardefinable1}(i), $L(\mu)=L\big({\mu_n}^{n\star}\big),$ hence $\mu\approx {\mu_n}^{n\star}$ by the Loeb measure construction.

Therefore, it holds for each $n\in\mathbb{N}$ that $\forall X\in{^*\sigma}\mathcal{B}\big( \vert\mu(X)-{\mu_n}^{n\star}(X)\vert<n^{-1}\big).$

By $\aleph_1$-saturation, extend $\{\mu_n\}_{n\in\mathbb{N}}$ to an internal sequence. Then it holds for all small enough $N\in{^*\mathbb{N}}$ that $\forall X\in{^*\sigma}\mathcal{B}\big( \vert\mu(X)-{\mu_{N}}^{N\star}(X)\vert<N^{-1}\big).$

Therefore, for any small enough infinite $N\in{^*\mathbb{N}},$ if we let $\nu=\mu_N,$ then $\mu \approx\nu^{N\star}.$

(ii) $\Rightarrow$ (iii): By $\aleph_1$-saturation, there are arbitrarily small infinite factorials.

(iii) $\Rightarrow$ (i): Suppose $\mu \approx\nu^{N!\star},$ where $N\in{^*\mathbb{N}}$ is infinite and $\nu\in{^*\mathcal{D}(\mathfrak{A})}.$ For each $n\in\mathbb{N},$ let $\mu_n=\nu^{N!/n\star}.$ Then $\mu\approx\mu_n^{n\star},$ hence $L(\mu)=L\big(\mu_n^{n\star}\big)=L(\mu_n)^{n\star}$ by Def.~\ref{neardefinable}(iv) and Prop.~\ref{neardefinable1}(i) again.
\end{proof}

Now we define a L\'{e}vy process along a timeline.

\begin{definition}\label{levyproc}
By a timeline, we mean an interval $I$ with endpoints from a linearly ordered semigroup $(S, +, \leq )$ such that the left endpoint of $I$ is denoted by $0$ and the right one by $1.$

(Note that the symbols $+$ and $0$ are used both for this linearly ordered semigroup and for the semigroup of $\mathfrak{A},$ but the intended meaning should be clear from the context.)

Let $\mu\in\mathcal{D}(\mathfrak{A})$ be infinitely divisible. By a L\'{e}vy process corresponding to $\mu$ \emph{w.r.t.} $I$ we mean a mapping $X: I \longrightarrow\mathcal{D}(\mathfrak{A})$ such that
\[X(0)=\delta_0, X(1)=\mu\quad\text{and}\quad X(s+t)= X(s)\star X(t)\quad\text{for all}\quad s, t, s+t \in I.\]

A L\'{e}vy process corresponding an infinitely divisible $\mu\in\Gamma({^*\mathfrak{A}})$ \emph{w.r.t.} $I$ is defined similarly.\nolinebreak\hfill $\Box$
\end{definition}

The L\'{e}vy process above can be regarded as an evolution along a ``straight line segment of probabilities" from the deterministic element $\delta_0$ to the random element $\mu.$ Intuitively, one expects that in general unless the geometry is complicated, there should be only one unique ``straight line segment" joining $\delta_0$ and $\mu.$ Of course this depends on the uniqueness of the $n^{\rm th}$ root of the convolution product.

Main examples of the $I$ considered are the real interval $[0,1]$ from $ \mathfrak{R},$ or $[0,1]\cap\mathbb{Q}$ from $(\mathbb{Q}, +, \cdot,\leq, 0, 1)$ or the hyperfinite timeline of the form $\big\{ 0,N^{-1}, \ 2N^{-1},\cdots, NN^{-1}=1\big\},$ identifiable with $\{0, 1, 2,\cdots, N\},$ from $({^*\mathbb{N}}, +, \cdot,\leq, 0, 1)$ for some infinite $N\in^*\mathbb{N}.$

\begin{remark}\label{levyinfdiv}
Given a L\'{e}vy process \emph{w.r.t.} $I$ of the above types, $X(t)$ is infinitely divisible for every $t\in I,$ since for any $n\in\mathbb{N}^+,$ by repeated applications of the additive condition to $X(n^{-1}t),$ we have $X(n^{-1}t)^{n\star}=X(t).$
\nolinebreak\hfill $\Box$
\end{remark}

If $X$ is a L\'{e}vy process \emph{w.r.t.} $[0,1],$ the restriction of $X$ to $[0,1]\cap\mathbb{Q}$ is of course also a L\'{e}vy process. But converting a L\'{e}vy process \emph{w.r.t.} an infinite hyperfinite timeline to a L\'{e}vy process \emph{w.r.t.} $[0,1]$ would require some continuity conditions. However it is simple to get a L\'{e}vy process in $\Gamma({^*\mathfrak{A}})$ \emph{w.r.t.} a hyperfinite timeline.

\begin{proposition}
Let $\mu\in\Gamma({^*\mathfrak{A}})$ be infinitely divisible.

{\rm (i)} For some infinite $N\in{^*\mathbb{N}}$ there exists a L\'{e}vy process corresponding to $\mu$ w.r.t. the hyperfinite timeline $I=\big\{nN^{-1}\mid  n=0,1,\cdots,N\big\}.$

{\rm (ii)} There exists a L\'{e}vy process corresponding to $\mu$ w.r.t. $I =\mathbb{Q}\cap [0,1].$
\end{proposition}

\begin{proof}
(i): By Prop.~\ref{infdiv1}(ii), $\mu=L\big(\nu^{N\star}\big)$ for some infinite $N\in{^*\mathbb{N}}$ and $\nu\in{^*\mathcal{D}(\mathfrak{A})}.$ Then we simply define $X: I\longrightarrow\Gamma({^*\mathfrak{A}})$ by $X(nN^{-1}) =L\big(\nu^{n\star}\big), n= 0, 1,\cdots, N.$

(ii): By Prop.~\ref{infdiv1}(iii), $\mu=L\big(\nu^{N!\star}\big)$ for some infinite $N\in{^*\mathbb{N}}$ and $\nu\in{^*\mathcal{D}(\mathfrak{A})}.$ Then we define the L\'{e}vy process $X: I\longrightarrow\Gamma({^*\mathfrak{A}})$ by letting $X(nm^{-1})=L\big(\nu^{nm^{-1}N!\star}\big),$ where $n,m\in\mathbb{N}$ with $nm^{-1}\in \mathbb{Q}\cap [0,1].$
\end{proof}

\begin{remark}\label{levyindependent}
In the case $\mathfrak{A}=\mathfrak{R}$ one can show for example by \cite{Ng2} that the definition of the L\'{e}vy processes above does not depend on a particular choice of the $\nu$ in the above proof. But we don't know how to identify general $\mathfrak{A}$ that satisfies this property. \nolinebreak\hfill $\Box$
\end{remark}

The more difficult problem is to find L\'{e}vy processes in $\mathcal{D}(\mathfrak{A})$ \emph{w.r.t.} the continuous timeline $I=[0,1].$

\begin{proposition}\label{infdiv2}
Let $\mu\in\mathcal{D}(\mathfrak{A})$ and $r\in{^*\mathbb{R}}$ have standard part $s\geq 0.$ Then $^\circ\big(e^{r({\s\mu}\star -1)}\big)=e^{s(\mu\star -1)}.$

In particular, $^\circ\big(e^{r({\s\mu}\star -1)}\big)\in\mathcal{D}(\mathfrak{A}).$
\end{proposition}

\begin{proof}
By $r$ being finite, $e^{r({\s\mu}\star -1)}\approx e^{-r}\sum_{n=0}^K\frac{r^n}{n!}{\s\mu}^{n\star}\approx e^{-s}\sum_{n=0}^K\frac{s^n}{n!}{\s\mu}^{n\star}$ for any infinite $K\in{^*\mathbb{N}}.$

Hence $L\big(e^{r({\s\mu}\star -1)}\big)$ coincide with $e^{s(\mu\star -1)}$ on $\sigma\mathcal{B},$ \emph{i.e.} $^\circ\big(e^{r({\s\mu}\star -1)}\big)=e^{s(\mu\star -1)}.$
\end{proof}

We need the following little fact before proving the next theorem.

\begin{proposition}\label{compare}
Let $\mu,\nu\in{^*\mathcal{D}(\mathfrak{A})}.$ Then $\mu\lessapprox\nu$ implies $\mu \approx\nu.$
\end{proposition}

\begin{proof}
Suppose $\mu\lessapprox\nu$ and there is $S\in{^*\sigma}\mathcal{B}$ such that $\mu(S)\lnapprox\nu(S).$ Consider the complement $S^\text{c},$ then $\mu (S^\text{c})=1-\mu(S)\gnapprox 1-\nu (S)=\nu(S^\text{c}),$ a contradiction.
\end{proof}

\begin{theorem}\label{exp}
Let $\mu\in{^*\mathcal{D}(\mathfrak{A})}$ and let $r\in{^*[0,\infty)}.$ Then for all large enough $K\in{^*\mathbb{N}}$ there is $\lambda\in{^*\mathcal{D}(\mathfrak{A})}$ such that $e^{r(\mu\star - 1)}\approx \lambda^{K\star}.$ Moreover, we can take
\[\lambda =\big(1+rK^{-1}\big)^{-1}\delta_0 + rK^{-1}\big(1+rK^{-1}\big)^{-1}\mu,\quad i.e. \quad
\big(1+rK^{-1}\big)^{-1}\big(\delta_0+rK^{-1}\mu\big).\]
\end{theorem}

\begin{proof}
Let $K\in{^*\mathbb{N}}$ and define $\lambda$ as above. By transferring Prop.~\ref{lincon}, as a convex combination, $\lambda\in{^*\mathcal{D}(\mathfrak{A})}.$ Then by transferring Prop.~\ref{combination},
\[\lambda^{K\star} = \big(1+rK^{-1}\big)^{-K} \sum_{n=0}^K{K \choose n}\frac{r^n}{K^n}\mu^{n\star}.\]
Note that
\[{K \choose n}\frac{1}{K^n} = \frac{K!}{K^n (K-n)!}\frac{1}{n!} = \prod_{i=0}^{n-1}\bigg(1-\frac{i}{K}\bigg)
 \frac{1}{n!} \leq  \frac{1}{n!}.\]
Note also that for all large enough $K\in{^*\mathbb{N}}$ we have $e^r \big(1+rK^{-1}\big)^{-K}\approx 1.$

Moreover, for such $K$ we have
\[e^{-r}\sum_{n=0}^K\frac{r^n}{n!} \mu^{n\star} \approx e^{r(\mu\star - 1)}.\]
Hence it follows that for large enough $K\in{^*\mathbb{N}}$ we have
\[\lambda^{K\star} \lessapprox e^{-r}\sum_{n=0}^K\frac{r^n}{n!} \mu^{n\star} \approx e^{r(\mu\star - 1)},\]
the result now follows from Proposition \ref{compare}.
\end{proof}

By Prop.~\ref{infdiv1}(iii), Thm.~\ref{exp} and Prop.~\ref{neardefinable1}(ii), we immediately have the following.

\begin{corollary}\label{infdiv3}
Let $\mu\in{^*\mathcal{D}(\mathfrak{A})}$ and let $r\in{^*[0,\infty)}.$ Then $e^{r(L(\mu)\star - 1)}$ is
infinitely divisible in $\Gamma({^*\mathfrak{A}}).$\nolinebreak\hfill$\Box$
\end{corollary}

\begin{remark}
In the proof of Thm.~\ref{exp}, the size of the $K$ depends on the $r$ only. So, if the $\lambda^{K\star}$ is given first, one can identify some $r$ and obtain $\mu$ so that $e^{r(\mu\star - 1)}\approx \lambda^{K\star}.$ Therefore this would have shown that in $\Gamma({^*\mathfrak{A}})$ infinite divisibility is equivalent to being a convolution exponential. However, such $\mu$ would be a (signed) real-valued measure instead of a probability measure. Perhaps, this is enough justification to extend the framework here to more general measures.\nolinebreak\hfill$\Box$
\end{remark}

The following gives more information about convolution exponentials in $\Gamma({^*\mathfrak{A}}).$

\begin{corollary}\label{expconv}
{\rm (i)} Let $\mu\in{^*\mathcal{D}(\mathfrak{A})}$ and $r,s\in{^*[0,\infty)}.$ Then
\[e^{r(L(\mu)\star -1)}\star e^{s(L(\mu)\star -1)}= e^{(r+s)(L(\mu)\star -1)}.\]

{\rm (ii)} Convolution exponentials in $\Gamma({^*\mathfrak{A}})$ are always infinitely divisible and have $n^{\rm th}$ roots given by convolution exponentials: Let $\mu\in{^*\mathcal{D}(\mathfrak{A})},$ $r\in{^*[0,\infty)}$ and $n\in\mathbb{N}^+.$ Then
\[\Big(e^{rn^{-1}(L(\mu)\star -1)}\Big)^{n\star} = e^{r(L(\mu)\star -1)}.\]
\end{corollary}

\begin{proof}
(i): Apply Thm.~\ref{exp} to $e^{(r+s)(\mu\star -1)}$ and choose large enough $K\in{^*\mathbb{N}}$ such that $\varepsilon = rsK^{-1}\approx 0.$ Then
\begin{align*}
e^{r (\mu\star -1)}\star e^{s (\mu\star -1)} &\approx
\bigg(1+\frac{r}{K}\bigg)^{-K}\bigg(\delta_0+\frac{r}{K}\mu\bigg)^{K\star} \star
\bigg(1+\frac{s}{K}\bigg)^{-K}\bigg(\delta_0+\frac{s}{K}\mu\bigg)^{K\star}\\
&= \bigg(1+\frac{r+s+\varepsilon}{K}\bigg)^{-K}\bigg(\delta_0+\frac{r+s}{K}\mu+\frac{\varepsilon}{K}\mu^{2\star}
\bigg)^{K\star}\\
&\approx \bigg(1+\frac{r+s}{K}\bigg)^{-K}\bigg(\delta_0+\frac{r+s}{K}\mu\bigg)^{K\star}\approx e^{(r+s)(\mu\star -1)}.
\end{align*}
The arrival of the first term following the second last $\approx$ is justified by $\varepsilon\approx 0$ and $\big(1+\frac{r+s}{K}\big)^{K}\approx e^{r+s},$ while the second term is justified by the following estimate: Let $X\in{^*\sigma\mathcal{B}},$ then
\begin{align*}
&\bigg\vert \bigg(\delta_0+\frac{r+s}{K}\mu+\frac{\varepsilon}{K}\mu^{2\star}\bigg)^{K\star}\big(X\big)
-\bigg(\delta_0+\frac{r+s}{K}\mu\bigg)^{K\star}\big(X\big)\bigg\vert\\
= &\bigg(\sum_{n=1}^K{K \choose n}\frac{\varepsilon^n}{K^n}\Big(\delta_0+\frac{r+s}{K}\mu\Big)^{(K-n)\star}\mu^{2n\star}\bigg)\big(X\big)
\leq \bigg(1+ \frac{\varepsilon}{K}\bigg)^K-1\approx e^\varepsilon-1\approx 0.
\end{align*}
Now by Prop.~\ref{neardefinable1}(i),(ii), we have
\begin{align*}
e^{r (L(\mu)\star -1)}\star e^{s (L(\mu)\star -1)} &=L\big(e^{r (\mu\star -1)}\big)\star L\big(e^{s (\mu\star -1)}\big)=L\big(e^{r (\mu\star -1)}\star e^{s (\mu\star -1)}\big)\\
          &=L\big(e^{(r+s)(\mu\star -1)}\big)=e^{(r+s)(L(\mu)\star -1)}.
\end{align*}

(ii): By first applying (i) with both $r,s$ replaced by $rn^{-1},$ then applying (i) again with $r,s$ replaced by $rn^{-1}, 2rn^{-1},$ ... ..., after iterating $(n-1)$ times, the conclusion follows.
\end{proof}

Now an almost identical result for convolution exponentials in $\mathcal{D}(\mathfrak{A}):$

\begin{corollary}\label{expconv0}
{\rm (i)} Let $\mu\in\mathcal{D}(\mathfrak{A})$ and $r,s\in [0,\infty).$ Then
\[e^{r(\mu\star -1)}\star e^{s(\mu\star -1)}= e^{(r+s)(\mu\star -1)}.\]

{\rm (ii)} Convolution exponentials in $\mathcal{D}(\mathfrak{A})$ are always infinitely divisible and have $n^{\rm th}$ roots given by convolution exponentials: Let $\mu\in\mathcal{D}(\mathfrak{A}),$ $r\in [0,\infty)$ and $n\in\mathbb{N}^n.$ Then
\[\Big(e^{rn^{-1}(\mu\star -1)}\Big)^{n\star} = e^{r(\mu\star -1)}.\]
\end{corollary}

\begin{proof}
(i): Apply the same computation in the proof of Cor.~\ref{expconv}(i) for ${^*\mu},$ we have
\[e^{r ({^*\mu}\star -1)}\star e^{s ({^*\mu}\star -1)}\approx e^{(r+s)({^*\mu}\star -1)}.\]
Let $X\in\sigma\mathcal{B},$ then, by transfer,
\[\big(e^{r(\mu\star -1)}\star e^{s(\mu\star -1)}\big)(X)=\big(e^{r({^*\mu}\star -1)}\star e^{s({^*\mu}\star -1)}\big)({^*X})\approx e^{(r+s)({^*\mu}\star -1)}({^*X})=e^{(r+s)(\mu\star -1)}(X).\]
Hence $\big(e^{r(\mu\star -1)}\star e^{s(\mu\star -1)}\big)(X)=e^{(r+s)(\mu\star -1)}(X)$ for all
$X\in\sigma\mathcal{B},$ \emph{i.e.} $e^{r(\mu\star -1)}\star e^{s(\mu\star -1)}= e^{(r+s)(\mu\star -1)}.$

(ii) follows by iterating (i), as in the proof of Cor.~\ref{expconv}(ii).
\end{proof}

The following gives the ``mode" of continuity in the convolution powers.

\begin{corollary}\label{expconv1}
Let $\lambda$ and $K$ as in Thm.~\ref{exp}. Let $L\in{^*\mathbb{N}}$ such that $rLK^{-1}\approx 0.$ Then $\lambda^{(K+L)\star}\approx\lambda^{K\star}.$
\end{corollary}

\begin{proof}
For notational convenience, we extend our definitions slightly and use real-valued measures. First note that
\[\lambda^{(K+L)\star}-\lambda^{K\star}=\lambda^{K\star}\star(\lambda^{L\star}-\delta_0)=\lambda^{K\star}\star (\lambda-\delta_0)\star\sum_{n=0}^{L-1}\lambda^{n\star}.\]
Since $\lambda-\delta_0=\big(1+rK^{-1}\big)^{-1}rK^{-1}(\mu-\delta_0),$ we have $\forall X\in{^*\sigma}\mathcal{B}\Big(\vert\lambda^{(K+L)\star}(X)-\lambda^{K\star}(X)\vert\lessapprox 2rLK^{-1}\approx 0\Big).$ \emph{i.e.} $\lambda^{(K+L)\star}\approx\lambda^{K\star}.$
\end{proof}

Now we are ready to show that, in both the case of $\Gamma({^*\mathfrak{A}})$ and the case of $\mathcal{D}(\mathfrak{A}),$ L\'{e}vy process \emph{w.r.t.} the continuous timeline $I=[0,1]$ always exists for a convolution exponential. The question of uniqueness is still open for the general case other than $\mathfrak{R}.$

\begin{theorem}\label{contlevy}
Let $I=[0,1],$ $[0,1]\cap\mathbb{Q}$ or any hyperfinite timeline of the form $\big\{0,N^{-1}, 2N^{-1},\cdots, 1\big\},$ where $N\in{^*\mathbb{N}}.$

{\rm (i)} Let $\mu\in\Gamma({^*\mathfrak{A}})$ be a convolution exponential. Then there exists a L\'{e}vy process $X:I\longrightarrow\Gamma({^*\mathfrak{A}})$ such that $X(1)=\mu.$

{\rm (ii)} Let $\mu\in\mathcal{D}(\mathfrak{A})$ be a convolution exponential. Then there exists a L\'{e}vy process $X:I\longrightarrow\mathcal{D}(\mathfrak{A})$ such that $X(1)=\mu.$
\end{theorem}

\begin{proof}
(i): For some $\nu\in{^*\mathcal{D}(\mathfrak{A})}$ and $r\in [0,\infty),$ we have $\mu=e^{r(L(\nu)\star-1)}.$

Then, for $I=[0,1],$ define $X:I\longrightarrow \Gamma({^*\mathfrak{A}})$ by $X(t)=e^{tr(L(\nu)\star-1)}$ and so, by Cor.~\ref{expconv}(i), $X$ is a L\'{e}vy process.

The restriction of $X$ to $[0,1]\cap\mathbb{Q}$ is clearly also a L\'{e}vy process.

For $I=\big\{0,N^{-1}, 2N^{-1},\cdots, 1\big\},$ we define $X(nN^{-1})$ as $L\big(e^{nN^{-1}r(\nu\star-1)}\big)$ and use Cor.~\ref{expconv}(i) again.

(ii): We write $\mu=e^{r(\nu\star-1)}$ for some $\nu\in\mathcal{D}(\mathfrak{A})$ and $r\in [0,\infty).$
Then define $X$ as in the proof of (i) but apply Cor.~\ref{expconv0}(i).

In the case $I=\big\{0,N^{-1}, 2N^{-1},\cdots, 1\big\},$ we define $X(nN^{-1})$ as $^\circ\big(e^{nN^{-1}r({^*\nu}\star-1)}\big),$ which is in $\mathcal{D}(\mathfrak{A}),$ by Prop.~\ref{infdiv2}.
\end{proof}

Now we formulate a weak converse of Cor.~\ref{expconv0}(ii) as an important property which basically says that infinitely divisible probabilities are in the closure of exponential ones in a certain sense. This property holds for the case of $\mathfrak{R}.$ Combined with Fourier analysis, the celebrated L\'{e}vy-Khintchine formula is a corollary to this property.

\begin{definition}\label{levykhintchine}
We say that $\mathfrak{A}$ has the L\'{e}vy-Khintchine property if for every $\mu\in\mathcal{D}(\mathfrak{A}),$ $\mu$ is infinitely divisible iff $\mu={^\circ\big(e^{r(\nu\star -1)}\big)}$ for some $\nu\in{^*\mathcal{D}(\mathfrak{A})}$ and $r\in{^*[0,\infty)}.$ \nolinebreak\hfill $\Box$
\end{definition}

\begin{remark}\label{lknote}
Suppose $\mathfrak{A}$ has the L\'{e}vy-Khintchine property. Then for each infinitely divisible $\mu\in\mathcal{D}(\mathfrak{A}),$ by Thm.~\ref{contlevy}(i), there is a L\'{e}vy process $X$ in $\Gamma({^*\mathfrak{A}})$ so that the restriction of $X(1)$ to $\sigma\mathcal{B}$ is $\mu.$ However, $X(1)$ needs not be $L({^*\mu})$ in general.\nolinebreak\hfill $\Box$
\end{remark}

We isolate the following property for an infinitely divisible probability which requires the $n^{\rm th}$ roots to be concentrate sufficiently near $0.$

\begin{definition}\label{concen}
We say that $\mu\in\mathcal{D}(\mathfrak{A})$ has the concentration property if there exists $\lambda\in{^*\mathcal{D}(\mathfrak{A})},$ $ r\in{^*[0,\infty)}$ and infinite $K\in{^*\mathbb{N}}$ such that

(i) $e^r\big(1+rK^{-1}\big)^{-K}\approx 1;$

(ii) $\forall X\in{^*\sigma}\mathcal{B}\big(\vert{^*\mu}(X)-\lambda^{K\star}(X)\vert\leq K^{-1}\big);$

(iii) $\lambda\left(\{0\}\right)\geq \big(1+rK^{-1}\big)^{-1}.$ \nolinebreak\hfill $\Box$
\end{definition}

Note that by Prop.~\ref{infdiv1}, for the above $\mu,$ $L({^*\mu})$ has to be infinitely divisible in $\Gamma({^*\mathfrak{A}}).$

From classical results such as those in \cite{Stroock} one can show that every infinitely divisible probability on
$\mathfrak{R}$ has the concentration property.

Our main interest of the property is the following:

\begin{theorem}
Suppose that every infinitely divisible probability in $\mathcal{D}(\mathfrak{A})$ has the concentration property.
Then $\mathfrak{A}$ has the L\'{e}vy-Khintchine property.
\end{theorem}

\begin{proof}
Let $\mu\in\mathcal{D}(\mathfrak{A})$ be infinitely divisible, with the $\lambda, K$ and $r$ given as in Def.~\ref{concen}. For some internal $\mathfrak{C}\succ\mathfrak{A},$ each $X\in{^*\mathcal{B}}$ has the form $\varphi^{C}(\cdot)$ for some ${^*{\mathcal{L}}}_{A}$-formula $\varphi.$ To define an internal probability $\nu$ on $^*\mathfrak{A},$ for $X=\varphi^{C}(\cdot),$ let
\[\nu(X) =
\left\{%
\begin{array}{ll}
    \big(1+Kr^{-1}\big)\lambda (X) - Kr^{-1}, &
    \hbox{if ${^*\mathfrak{A}}\models\varphi (0),$ \emph{i.e.} if $0\in X$} \\
     \\
    \big(1+Kr^{-1}\big)\lambda (X), & \hbox{otherwise}. \\
\end{array}%
\right.\]
Using the lower bound for $\lambda\left(\{0\}\right),$ it is
easy to check that $\nu$ extends uniquely to an internal probability on $^*\mathfrak{A}$ and belongs to ${^*\mathcal{D}(\mathfrak{A})}.$ We still call the extension as $\nu.$ Now we
can re-write $\lambda$ as the convex combination
\[\lambda = \big(1+rK^{-1}\big)^{-1} \big(\delta_0+rK^{-1}\nu\big).\]
Then similar to the proof of Thm.~\ref{exp}
\[{^*\mu}\approx\lambda^{K\star}\lessapprox e^{-r}\sum_{m=0}^K\frac{r^m}{m!} \nu^{m\star}\approx e^{r(\nu\star - 1)},\]
and hence ${^*\mu}\approx e^{r(\nu\star -1)}$ and so $\mu={^\circ{L({^*\mu})}}={^\circ\big(e^{r(\nu\star -1)}\big)}.$
\end{proof}

\begin{remark}\label{conjecture}
{(i)} We conjecture that $p$-adic fields have the L\'{e}vy-Khintchine property.

{(ii)} In our framework, we do not exclude the possibility that the semigroup from $\mathfrak{A}$ could be compatible with the semigroup that indexes a L\'{e}vy process. It may be worthwhile to formulate L\'{e}vy processes which are definable processes in some sense.\nolinebreak\hfill $\Box$
\end{remark}

\end{document}